\documentclass[11pt]{article}
\usepackage{graphicx} 
\usepackage{amsmath, amsthm, amssymb}
\usepackage{mathtools}
\usepackage[parfill]{parskip}
\usepackage{soul}
\usepackage{multirow}
\usepackage{caption}
\usepackage{subcaption}
\usepackage{mathtools}
\usepackage{hyperref}
\usepackage{epsfig,fullpage,psfrag,url,verbatim}
\usepackage{graphics}
\usepackage{graphicx}
\usepackage{epstopdf}
\usepackage{amsfonts}
\usepackage{array}
\usepackage{booktabs}
\usepackage{tikz}
\usepackage{tikz-cd}
\usepackage{pgf}
\usepackage{rotating}
\usepackage{tikz}
\usepackage{bbm}
\usepackage{longtable}
\usepackage{authblk}

\newtheorem{theorem}{Theorem}%
\newtheorem{proposition}[theorem]{Proposition}%

\newtheorem{lemma}{Lemma}
\newtheorem{corollary}{Corollary}

\allowdisplaybreaks

\usepackage{mathrsfs}
\usepackage{indentfirst}
\usepackage{float}
\usepackage{multirow}
\usepackage[linesnumbered,ruled]{algorithm2e}
\usepackage{algpseudocode} 
\usepackage{diagbox}

\begin{document}

\title{Modeling and solving cascading failures across interdependent infrastructure systems}
\author[1]{Yijiang Li}
\author[1]{Kibaek Kim}
\author[1]{Sven Leyffer}
\author[1]{Matt Menickelly}
\author[2]{Lawrence Paul Lewis}
\author[2]{Joshua Bergerson}
\affil[1]{Mathematics and Computer Science, Argonne National Laboratory}
\affil[2]{Decision and Infrastructure Sciences, Argonne National Laboratory}
\affil[ ]{\{yijiang.li, kimk, leyffer, mmenickelly, plewis, jbergerson\}@anl.gov}
\date{}

\maketitle

\begin{abstract}
Physical infrastructure systems supply crucial resources to residential, commercial, and industrial activities. These infrastructure systems generally consist of multiple types of infrastructure assets that are interdependent. In the event of a disaster, some of the infrastructure assets can be damaged and disabled, creating failures that propagate to other assets that depend on the disabled assets and cause a cascade of failures that may lead to a potential system collapse. We present a bilevel interdiction model in this paper 
to study this problem of cascading failures in a system of interdependent infrastructure systems with a nondeterministic dependency graph. We also propose a computationally tractable reformulation of the proposed bilevel model and utilize a Benders-type decomposition algorithm to solve the resulting formulation. Computational experiments are performed using infrastructure networks generated from anonymized real-world data to validate the performance of this algorithm. 
\end{abstract}

\maketitle

\section{Introduction}\label{sec:intro}

Physical infrastructure systems provide crucial resources such as power, water, and internet services for residential, commercial, and industrial activities. An infrastructure system generally consists of multiple types of assets and service connections. These infrastructure assets can be damaged and disabled in the event of a disaster. For example, a hurricane could leave particular infrastructure assets inoperable (i.e., disabled). 
In addition, physical infrastructure systems are generally interdependent. A cascading failure occurs when one inoperable asset induces the failure of another asset that depends on it. Such failures, if not handled properly, could propagate further and lead to an overall system collapse. These cascading failures may significantly increase the geographic footprint of infrastructure failures and overall losses beyond the initial damages resulting from an initiating event (e.g., hurricane). Cascading failures are discussed in \cite{PhysRevE.66.065102, Chang2009, Leavitt2006}. 

This situation motivates us to study a system of interdependent and interconnected infrastructure systems and identify the most severe contingencies, or combination of assets, whose failures result in the most severe impact on the service level of the network. This assessment illustrates how cascading failures affect the overall system and potentially scale up the consequences of an initial disruption event. The identified set of infrastructure assets is essential to the network's functionality and deserves special attention during blue sky conditions and following the occurrence of natural disasters. 

We propose a network model that does not assume a given description of the connectivity between the infrastructure assets. In particular, we assume given and fixed dependencies between asset classes as summarized in Table~\ref{table:asset_depedency}. 
The table presents data that was anonymized from real-world data collected from infrastructure owners and operators by the Decision and Infrastructure Sciences (DIS) division at Argonne National Laboratory in support of recovery efforts in Puerto Rico following Hurricane Maria. These infrastructure systems consist of multiple assets from each  asset class, with directed edges between assets based on dependency relationships between asset classes. This configuration allows us to capture any dependencies as natural consequences of asset classes while providing flexibility in modeling the dependencies beyond the asset classes based on additional factors that can be potentially unknown or stochastic in nature. 

\begin{table}[ht]
\caption{Asset information} \label{table:asset_depedency}
	\centering
 \small 
	\begin{tabular}{|c|c|}
	\hline
	Asset class & Upstream dependencies \\\hline\hline
	Cell tower & Substation \\\hline
	EMS & Substation, Cell tower, Water, Health care\\\hline
	Thermal generation & Transport, Cell tower \\\hline
	Renewable generation & Cell tower\\\hline
	Substation & Thermal or Renewable generation, Cell tower\\\hline
	Health care & Transport, EMS, Substation, Cell tower, Water\\\hline
	Pharmacies & Transport, Substation, Cell tower, Water\\\hline
	Transport & Cell tower, Substation, Water\\\hline
	Water & Substation, Transport, Cell tower \\\hline
	\end{tabular}
\end{table}

We pose the problem of identifying the most severe contingencies as a bilevel interdiction problem in the attacked-defender (AD) framework (see a general overview in \cite{Wood2011}) where the attacker (i.e., natural disaster), also the leader in the bilevel problem, seeks to disrupt a set of assets that maximizes the disruption to the network.
The defender (i.e., network operator), also the follower in the bilevel problem, maximizes the remaining service level of the network given the attacked assets and cascade of failures.

The rest of this paper is organized as follows. Section \ref{sec:lit_review} provides a literature review of relevant work. Section \ref{sec:model} presents the bilevel interdiction model, while Section \ref{sec:bender} demonstrates the decomposition algorithm. Section \ref{sec:numerical} presents the numerical experiments.  Section \ref{sec:conclusion} summarizes our work and briefly suggests avenues for future research. 

\section{Literature review} \label{sec:lit_review}

The study of cascading failures in networks is a key part of the study of resilient networks. 
For a detailed review of the resilience of networks, we refer to \cite{Mishra_2020}. 
We first give a summary of the various ways to model this problem. 
For this presentation, we broadly divide 
relevant work into two sets: a set of work that uses tools that are not optimization-based and a set of tools that are optimization-based. 

The first set of work focuses on an abstract network model and studies the effect of cascading failures across the network via tools not using optimization. 
\cite{Duan_2019} study the propagation of cascading failures in networks of different characteristics and present a study accounting for the dynamics of the nodes in the network and analyzing the impact of various degrees of disruption to the nodes in the network. 
\cite{Smolyak_2020} consider a similar model for the cascading condition and propose an algorithm to select a set of more instrumental nodes in the failure propagation process. 
Moreover, they consider additional parameters that affect cascading failures and run simulations to validate the proposed algorithm. 
\cite{Duenas-Osorio_2009} propose a new way to model the overloads due to cascading failures in a risk and reliability assessment of complex infrastructure systems. 
The effect of the proposed model is illustrated for power systems. 

The second set of work utilizes optimization tools to tackle the problems and usually focuses on electric power systems. This work includes \cite{Dandurand2021, Zhao_2013, Salmeron_2015, Salmeron_2009, Sundar_2021, Arroyo_2010, Seyyedi_2023, Wang_2014, Johnson2022}. Most of these studies use the Stackelberg game framework, a bilevel formulation, to model the problem on an electric power network. 
In general, in the bilevel formulation, the attacker (i.e., leader) decides which set of nodes to disable, and the defender (i.e., follower) decides on the operation of the network to minimize the network disruption of the given attack. 
In particular, \cite{Brown_2006} present the general attacker-defender (AD) framework that results in a bilevel formulation.
In addition, the bilevel AD approach can be extended into a trilevel formulation with fortification by the defender coming before the decisions of the attacker. 
Examples of such works include \cite{Alguacil2014} and \cite{mastin2023}. 
These authors consider applications in electric power systems and petroleum reserve systems, respectively.

Previous work has focused on exploiting particular features of the follower problem; this is especially true with optimal power flow (OPF) problems. 
In particular, \cite{Dandurand2021} consider a nonconvex alternating current optimal power flow (ACOPF) as its follower problem. 
ACOPF is inherently challenging to solve, and the authors propose a branch-and-bound framework with a convex semi-definite programming  relaxation to the follower problem to solve the bilevel problem. \cite{Salmeron_2015} consider recovery transformers in a power grid, which are replacement transformers that can be quickly set up for the attacked transformers, The authors study the value of these transformers in reducing network impacts of attacks. 
Another work \cite{Zhao_2013} incorporates an additional operation, line switching, in the follower problem. 
The authors perform a vulnerability analysis of the power grid and study the mitigation effect on line switching operations. 
\cite{Johnson2022} consider both the theoretical and empirical implications of relaxations and restrictions of the attacker-defender framework in modeling responses to attacks in transmission grids using direct current optimal power flow to speed up the solution time. 


To the best of our knowledge,  no previous studies have focused on quantifying the significance of network nodes to overall network performance in a general system of interdependent infrastructure systems, with explicit consideration of a stochastic/nondeterministic network topology. 

The remainder of this section provides a review of methods and algorithms for solving bilevel optimization problems. 
Because the attacker's decisions are typically modeled by binary variables, the resulting problem is a bilevel mixed-integer problem. 
We give a summary of relevant techniques and refer readers to  \cite{Kleinert_2021} for more details. 
In general, there are two types of algorithms for solving bilevel mixed-integer problems.
The first entails reformulating the bilevel problem into a single-level problem, while the second involves decomposition and the application of iterative algorithms. 

A common single-level reformulation is the value function reformulation and the associated \emph{high-point relaxation}, in which the reformulation maintains all the bilevel constraints but relaxes the optimality condition for the follower problem's objective function. 
Cuts are generated iteratively to cut off any solutions that are infeasible with respect to the optimality condition. 
Work by \cite{Fischetti_2020} proposes a new class of intersection cuts that is proven to be effective in cutting off infeasible points and reducing the size of the branch-and-bound tree generated by the high-point relaxation.
The authors validate the performance of the reformulation by considering common bilevel benchmark networks. 
Another single-level reformulation utilizes the Karush--Kuhn--Tucker  conditions to rewrite the follower problem.
\cite{Arroyo_2010} considers both such techniques to reformulate the bilevel problem into a single-level problem. 
This reformulation can be effective when the follower problem is relatively easy, for example, when it is a linear optimization problem. 

The second type of algorithm we review is a decomposition and iterative algorithm.
\cite{Tang_2016} construct an approximate leader problem from previously explored leader problem solutions and propose an iterative algorithm to solve the construction. 
When the follower problem satisfies certain assumptions, the proposed algorithm can be further improved. 
The authors demonstrate the effectiveness of the algorithm with min-max binary knapsack problems and min-max clique problems. 
\cite{Taninmis_2021} follow up on the work of \cite{Tang_2016} with further improvements to the latter's algorithm.
\cite{Johnson2022} propose a novel relaxation that leads to certain problem structures that scale well with problem size.
Benders decomposition algorithm was adapted by \cite{Salmeron_2009} where they compute a cut for each explored leader problem solution. 
The authors consider power or electric systems as applications of the proposed procedure. 
Another Benders decomposition algorithm was developed by \cite{israeli.wood:02} for the network-interdiction problem. 
The authors propose a class of cuts called super-valid inequalities that may cut off some feasible solutions but not the optimal solutions. 
In addition, others (e.g., \cite{Lozano_2022}) have investigated strategies for exploring the feasible solution space of the leader problem. 
Meta-heuristics have also been applied. 
For networks, a modified binary particle swarm optimization  was proposed by \cite{Seyyedi_2023} to solve the bilevel formulation in electric power systems. 

As for the extension of bilevel problems to trilevel problems to include fortification, \cite{Alderson2011} consider a decomposition algorithm for the general defender-attacker-defender model that optimizes decisions in an alternating fashion between the defender and attacker. 
The algorithm is demonstrated on a transportation problem derived from The Seven Bridges of K\"{o}nigsberg. \cite{Alguacil2014} propose an implicit enumeration scheme over fortification plans; once a fortification plan is determined and fixed, existing algorithms for bilevel optimization problems are applied. 
The success of such an enumeration scheme crucially depends on the size of the space of feasible fortification plans. 

\section{Problem formulation}\label{sec:model}

Our goal is to identify the most severe contingencies, or combinations of failures among assets, in terms of their impact on the service level of the network.
We formulate this problem as a bilevel interdiction problem, or a Stackelberg game. In our model, the leader (e.g., extreme weather condition) selects assets to disable so that the remaining service level of the network is minimized. The follower  (e.g., a network operator) seeks to maximize the service level subject to the disabled assets, taking into account the cascade of failures through the network. The follower problem is a multistage optimization problem, where each stage models a progression of the cascade. For simplicity we assume the progression of the cascade to be synchronized. 

\subsection{Notation} 
We use the convention that sets are denoted by calligraphic parameters; parameters and constants are denoted by upper-case characters; and optimization variables are denoted by lower-case characters.
\subsubsection{Parameters}
\begin{itemize}
    \item ${\cal F}$ is the set of assets (nodes).
    \item $N := |{\cal F}|$ is the total number of assets.
    \item $N_c \leq N$ is the maximum number of assets that can be disabled by the leader.
    \item $N_p \geq 1$ is the number of stages of failure cascades that we model in the follower problem.
    \item ${\cal A} \subset {\cal F} \times {\cal F}$ is the set of directed arcs that indicate dependencies; that is, $(s,f) \in {\cal A}$ indicates that the asset class of $f$ depends on the asset class of $s$ as given in Table~\ref{table:asset_depedency}.
    \item $W_f \geq 0$ is the importance weight of asset $f \in {\cal F}$.
    \item $P_{s,f,i}$ are the uncertainty parameters that represent the arc weights measuring the dependencies between assets $s$ and $f$ for $(s,f) \in {\cal A}$ in stage $i \in \{1,\ldots,N_p\}$, accounting for factors including the asset classes. We use $\boldsymbol{P_{f,i}}$ to denote the vector that contains the variables $P_{s,f,i}$. A generic polytope $\mathcal{P}_f$ of $P_{s,f,i}$ for an asset $f \in \mathcal{F}$ is given as
    \begin{equation}
    \mathcal{P}_f := \left\{P_{s,f,i} \ge 0: \sum_{s \in \mathcal{F}} U^f_{l,s}P_{s,f,i} \ge u^f_l, \; l \in L\right\}, \label{eq:P_f}
    \end{equation}
    where $U^f_{l,s}$ and $u^f_l$ are given constants. $L$ is an index set for all constraints in the polytope $\mathcal{P}_f$. In particular, we use the inequalities $\sum_{s \in \mathcal{F}} U^f_{l,s}P_{s,f,i} \ge u^f_l$ for $l \in L$ to account for dependencies as a result of the asset classes and other additional factors on $P_{s,f,i}$, so the description~\eqref{eq:P_f} remains generic and robust. Also note that we assume $\mathcal{P}_f$ is invariant across the stages of failure cascade. 
\end{itemize}
\subsubsection{Decision variables}
\begin{itemize}
    \item $x_f \in \{0,1\}$ are the upper-level decisions, where $x_f$ is equal to $1$ if asset $f$ is disabled and $0$ otherwise. We use $\boldsymbol{x}$ to denote the vector that contains all $x_f$. 
    \item $y_{f,i} \in [0,1]$ are the decision variables in stage $i$ for $i \in \{1,\ldots,N_p\}$. We interpret a variable $y_{f,i}$ as the fractional service level for asset $f$ in stage $i$. Each $y_{f,i}$ is jointly determined by the arc weights between $f$ and any assets on which it depends, as well as the remaining service levels of these assets from the previous stage. We use $\boldsymbol{y_{i}}$ to denote the vector that contains all $y_{f,i}$ for stage $i$. 
     
\end{itemize}



\subsection{Optimization model}
With this notation established, we can now define the leader and follower problems. We first consider the first stage of failure cascade. For a particular asset $f$, with arc weights $P_{s,f,1} \in \mathcal{P}_f$ and $x_s$ for $s \in \mathcal{F}$, we can compute the total loss of service level from the upstream dependencies as $\sum_{s: (s,f) \in \mathcal{A}}P_{s,f,1}x_s$. 
Next, assuming linear degradation of service level, we can model the service level via the constraint
\begin{equation}
	y_{f,1} \le \phi_f(\boldsymbol{P_{f,1}},\boldsymbol{x}) := \left[1- \sum_{s: (s,f) \in \mathcal{A}}P_{s,f,1}x_s\right]^+, \; \forall \boldsymbol{P_{f,1}} \in \mathcal{P}_f, \label{capacityreduction}
\end{equation}
where $[\cdot]^+ := \max\{\cdot,0\}$. Note that \eqref{capacityreduction} provides robustness to the optimization model on the factors affecting the interdependencies of assets that are unknown or stochastic. With this notation, we are ready to give our multistage formulation of the problem with a budget constraint for the leader:
\begin{subequations}
\label{eq:original}
\begin{align}
\min_{\boldsymbol{x} \in \{0,1\}^N} \; Q(\boldsymbol{x})
\quad \text{subject to} \quad e^\top \boldsymbol{x} \le N_c, \label{constr_budget}
\end{align}
where $e$ is the vector of all ones and  $Q(\boldsymbol{x})$ is given by the value of the follower problem: 
\begin{align}
Q(\boldsymbol{x}) := \max_{0 \le \boldsymbol{y} \le 1} \;\;\;\;& \sum_{f \in \mathcal{F}}\sum_{i=1}^{N_p}W_f y_{f,i}	\\
\text{subject to} \;\;\;\;& y_{f,1} \le 1 - x_f, \quad  f \in \mathcal{F} \label{inner_constr_capacity_first_direct}\\
& y_{f,i} \le y_{f,i-1}, \quad  f \in \mathcal{F}, i \in \{2,\ldots,N_p\} \label{inner_constr_capacity_multistage_direct}\\
& y_{f,1} \le \phi_f(\boldsymbol{P_{f,1}}, \boldsymbol{x}), \quad  f \in \mathcal{F}, \boldsymbol{P_{f,1}} \in \mathcal{P}_f \label{inner_constr_capacity_first}\\
& y_{f,i} \le \phi_f(\boldsymbol{P_{f,i}}, \boldsymbol{y_{i-1}}), \quad  f \in \mathcal{F}, i \in \{2,\ldots,N_p\}, \boldsymbol{P_{f,i}} \in \mathcal{P}_f, \label{inner_constr_capacity_second}
\end{align}
\end{subequations}
where the functions $\phi_f$ in the constraints~\eqref{inner_constr_capacity_first} and \eqref{inner_constr_capacity_second} are defined in expression~\eqref{capacityreduction}. Note that the proposed formulation \eqref{eq:original} is generic to model special cases such as deterministic arc weights (i.e., by singleton $\mathcal{P}_f$) and static arc weights (i.e., by $\boldsymbol{P}_{f,i}=\boldsymbol{P}_f$).

We observe that constraints \eqref{inner_constr_capacity_first} and \eqref{inner_constr_capacity_second} can be written as
\begin{subequations}
\label{min_constr}
\begin{align}
  & y_{f,1} \le \min_{\boldsymbol{P_{f,1}} \in \mathcal{P}_f} \phi_f(\boldsymbol{P_{f,1}},\boldsymbol{x}), \quad f \in \mathcal{F}, \label{capacityreduction_worstcase} \\
  & y_{f,i} \le \min_{\boldsymbol{P_{f,i}} \in \mathcal{P}_f} \phi_f(\boldsymbol{P_{f,i}}, \boldsymbol{y_{i-1}}), \quad  f \in \mathcal{F}, i \in \{2,\ldots,N_p\}. \label{inner_constr_capacity_second_worst_case}
\end{align}
\end{subequations}
This allows the following reformulation.
\begin{proposition}
    The follower problem $Q(\boldsymbol{x})$ can be equivalently written as follows:
    \begin{subequations}
    \label{nonlinear_follower}
    \begin{align}
        \max_{\boldsymbol{y},\boldsymbol{\lambda},\boldsymbol{v}} \quad & \sum_{f \in \mathcal{F}}\sum_{i=1}^{N_p}W_f y_{f,i}	\\
        \text{subject to} \quad 
        & \eqref{inner_constr_capacity_first_direct}-\eqref{inner_constr_capacity_multistage_direct}, \notag \\
        & y_{f,i} \le \boldsymbol{u^f}^\top \boldsymbol{v_{f,i}} + \lambda_{f,i}, \quad  f \in \mathcal{F}, i \in \{1,\ldots,N_p\} \label{nonlinear_follower:cons1}\\
        & \boldsymbol{U^f_{s}}^\top \boldsymbol{v_{f,1}} + \lambda_{f,1}x_s \le 0, \quad  f, s \in \mathcal{F}, (s,f) \in \mathcal{A} \label{nonlinear_follower:cons2}\\
        & \boldsymbol{U^f_{s}}^\top \boldsymbol{v_{f,i}} + \lambda_{f,i}(1-y_{s,i-1}) \le 0, \quad  f, s \in \mathcal{F}, (s,f) \in \mathcal{A}, i \in \{2,\ldots,N_p\}\\
        & \boldsymbol{U^f_{s}}^\top \boldsymbol{v_{f,i}} \le 0, \quad  f, s \in \mathcal{F}: (s, f) \notin \mathcal{A}, i \in \{1,\ldots,N_p\} \label{nonlinear_follower:cons4}\\
        & 0 \le \boldsymbol{v_{f,i}}, \quad f \in \mathcal{F}, i \in \{1,\ldots,N_p\} \label{nonlinear_follower:cons5}\\
        & 0 \le y_{f,i} \le 1, \quad f \in \mathcal{F}, i \in \{1,\ldots,N_p\},\label{nonlinear_follower:cons6}\\
        & 0 \le \lambda_{f,i} \le 1, \quad f \in \mathcal{F}, i \in \{1,\ldots,N_p\}, \label{nonlinear_follower:cons7}
    \end{align}
    \end{subequations}
    where $\boldsymbol{u^f}$, $\boldsymbol{v_{f,1}}$, and $\boldsymbol{U_s^f}$ of length $L$ are the vectors that contain the entries of $u^f_l$, $v_{f,1,l}$, and $U_{l,s}^f$ for $l \in L$, respectively.
\end{proposition}

\proof{Proof.}
First, the right-hand side of constraint~\eqref{capacityreduction_worstcase} can be equivalently written as the following problem for each $f\in\mathcal{F}$:
\begin{subequations}
\label{worst_case}
\begin{align}
\min_{\xi_{f,1},P_{s,f,1}} \quad & \xi_{f,1} \label{worst_case_obj}\\
\text{subject to} \quad & \xi_{f,1} \ge 1- \sum_{s: (s,f) \in \mathcal{A}}P_{s,f,1}x_s \;\;\; (\lambda_{f,1}) \label{worst_case_constr_1}\\
& \xi_{f,1} \ge 0 \label{worst_case_constr_2}\\
& P_{s,f,1} \ge 0 \label{worst_case_constr_3}\\
& \sum_{s \in \mathcal{F}} U^f_{l,s} P_{s,f,1} \ge u^f_l, \quad l \in L \;\;\; (v_{f,1,l}) \label{worst_case_constr_5},
\end{align}
\end{subequations}
where $\lambda_{f,1}$ and $v_{f,1,l}$ are the dual multipliers with respect to the corresponding constraints. The dual of \eqref{worst_case} is given by
\begin{subequations}
\begin{align}
\max_{\lambda_{f,1},\boldsymbol{v_{f,1}}} \quad & \boldsymbol{u^f}^\top \boldsymbol{v_{f,1}} + \lambda_{f,1}\\
\text{subject to} \quad & \boldsymbol{U^f_{s}}^\top \boldsymbol{v_{f,1}} + \lambda_{f,1}x_s \le 0, \quad  s \in \mathcal{F}, (s,f) \in \mathcal{A} \label{worst_case:constr_first_stage_dep}\\
& \boldsymbol{U^f_{s}}^\top \boldsymbol{v_{f,1}} \le 0, \quad  s \in \mathcal{F}: (s, f) \notin \mathcal{A} \label{worst_case:constr_first_stage_nondep}\\
& 0 \le \boldsymbol{v_{f,1}}, 0 \le \lambda_{f,1} \le 1. \label{worst_case:constr_first_stage_bound}
\end{align}
\end{subequations}
Consequently, we arrive at an equivalent reformulation of the constraint~\eqref{capacityreduction_worstcase} as follows (see, e.g., Section 3.3.3 in \cite{Nemirovski_lp}):
\begin{subequations}
\label{worst_case:constr_first_stage}
\begin{align}
& y_{f,1} \le \boldsymbol{u^f}^\top \boldsymbol{v_{f,1}} + \lambda_{f,1}, \quad  f \in \mathcal{F}\\
& \eqref{worst_case:constr_first_stage_dep} - \eqref{worst_case:constr_first_stage_bound}, \quad  f \in \mathcal{F}.
\end{align}
\end{subequations}

Similarly, rewriting the right-hand side of \eqref{inner_constr_capacity_second_worst_case}, we obtain the following equivalent formulation:
\begin{subequations}
\label{worst_case:constr_non_first_stage}
\begin{align}
& y_{f,i} \le  \boldsymbol{u^f}^\top \boldsymbol{v_{f,i}} + \lambda_{f,i}, \quad f \in \mathcal{F}, i \in \{2,\ldots,N_p\}\\
& \boldsymbol{U^f_{s}}^\top \boldsymbol{v_{f,i}} + \lambda_{f,i}(1-y_{s,i-1}) \le 0, \quad  f, s \in \mathcal{F}, (s,f) \in \mathcal{A}, i \in \{2,\ldots,N_p\} \label{worst_case:constr_non_first_stage_dep}\\
& \boldsymbol{U^f_{s}}^\top \boldsymbol{v_{f,i}} \le 0, \quad  f, s \in \mathcal{F}: (s, f) \notin \mathcal{A}, i \in \{2,\ldots,N_p\} \label{worst_case:constr_non_first_stage_nondep}\\
& 0 \le \boldsymbol{v_{f,i}}, 0 \le \lambda_{f,i} \le 1, \quad f \in \mathcal{F}, i \in \{2,\ldots,N_p\}. \label{worst_case:constr_non_first_stage_bound}
\end{align}
\end{subequations}
Therefore, constraints \eqref{min_constr} can be equivalently written by the set of constraints \eqref{worst_case:constr_first_stage} and \eqref{worst_case:constr_non_first_stage}.
\endproof

Note, however, that constraint~\eqref{worst_case:constr_non_first_stage_dep} is bilinear with respect to $\lambda$ and $y$, which can be linearized by using the McCormick envelopes with the substitution $\zeta_{f,s,i} = \lambda_{f,i} y_{s,i-1}$ for $i \in \{2,\ldots,N_p\}$. The McCormick envelopes result in the following set of constraints: 
\begin{subequations}
\label{mc}
\begin{align}
& \boldsymbol{U^f_{s}}^\top \boldsymbol{v_{f,i}} + \lambda_{f,i} - \zeta_{f,s,i} \le 0, \quad  f, s\in \mathcal{F}, (s,f) \in \mathcal{A}, i \in \{2,\ldots,N_p\} \label{mc_1}\\
& \zeta_{f,s,i} \ge 0, \quad  f, s\in \mathcal{F}, (s,f) \in \mathcal{A}, i \in \{2,\ldots,N_p\} \label{mc_2}\\
& \zeta_{f,s,i} \ge \lambda_{f,i} + y_{s,i-1} -1, \quad  f, s\in \mathcal{F}, (s,f) \in \mathcal{A}, i \in \{2,\ldots,N_p\} \label{mc_3}\\
& \zeta_{f,s,i} \le \lambda_{f,i}, \quad  f, s\in \mathcal{F}, (s,f) \in \mathcal{A}, i \in \{2,\ldots,N_p\} \label{mc_4}\\
& \zeta_{f,s,i} \le y_{s,i-1}, \quad  f, s \in \mathcal{F}, (s,f) \in \mathcal{A}, i \in \{2,\ldots,N_p\}. \label{mc_5}
\end{align} 
\end{subequations}

We make a remark concerning the McCormick envelopes~\eqref{mc}. 
Generally, McCormick envelopes lead to a relaxation because the variables $y_{f,i}$ and $\lambda_{f,i}$ are continuous between $0$ and $1$. In what follows, we establish the exactness of the McCormick envelopes at the optimal solutions of the follower problem. To that end, let $(\xi_{f,i}^*,P_{s,f,i}^*,\lambda_{f,i}^*,\boldsymbol{v}_{f,i}^*)$ be the optimal primal and dual solution, respectively, of problem \eqref{worst_case}. First, we prove that optimal dual solution $\lambda_{f,i}^*$ is binary. 

\begin{lemma} \label{lemma:lambda_value}
For $i \in \{1,\ldots,N_p\}$ and $f \in \mathcal{F}$, $\lambda_{f,i}^* \in \{0,1\}$.  
\end{lemma}

\proof{Proof of Lemma~\ref{lemma:lambda_value}.}
We consider $i = 1$; each case in which $i > 1$ will follow analogous arguments. For each $f \in \mathcal{F}$, the optimal value $\xi_{f,1}^*$ of problem \eqref{worst_case} is the larger value of $1- \sum_{s: (s,f) \in \mathcal{A}}P_{s,f,1}^* x_s$ and 0. The first-order optimality condition of problem \eqref{worst_case} with respect to $\xi$ is given by
\begin{equation}
    1 - \lambda_{f,1} - \pi_{f,1} = 0, \label{xi_first_order}
\end{equation}
where $\pi_{f,1} \ge 0$ represents the dual variable corresponding to the constraint~\eqref{worst_case_constr_2}. 
We denote its optimal value by $\pi_{f,1}^*$.
We consider both cases for $\xi_{f,1}^*$:
\begin{itemize}
\item When $\xi_{f,1}^* = 1- \sum_{s: (s,f) \in \mathcal{A}}P_{s,f,1}^*x_s$, we have $\xi_{f,1}^* > 0$. By complementary slackness, we have that $\pi_{f,1}^* = 0$ and \eqref{xi_first_order} gives $\lambda_{f,1}^*=1$.
\item When $\xi_{f,1}^* = 0$, we have $\xi_{f,1}^* > 1- \sum_{s: (s,f) \in \mathcal{A}}P_{s,f,1}^*x_s$; and thus by complementary slackness, we have $\lambda_{f,1}^* = 0$.
\end{itemize}
\endproof

Now we can prove the exactness of the McCormick envelopes at the optimal solutions of the follower problem. Formally, we have the following theorem. 
\begin{theorem} \label{theorem:mccormick_equiv}
    If we restrict $\lambda_{f,i} \in \{0,1\}$ for $i \in \{1,\ldots,N_p\}$ and $f \in \mathcal{F}$, the McCormick envelopes \eqref{mc_2}--\eqref{mc_5} are equivalent to the bilinear constraints $\zeta_{f,s,i} = \lambda_{f,i} y_{s,i-1}$ for $i \in \{2,\ldots,N_p\}$ and $f \in \mathcal{F}$. 
\end{theorem}
\proof{Proof of Theorem~\ref{theorem:mccormick_equiv}.}
Using Lemma~\ref{lemma:lambda_value}, we consider the two cases for each $\lambda_{f,i}$:
\begin{itemize}
    \item When $\lambda_{f,i} = 0$, $\zeta_{f,s,i} = \lambda_{f,i} y_{s,i-1}$ gives $\zeta_{f,s,i} = 0$ and $y_{s,i-1}\in[0,1]$ for $s \in \mathcal{F}$ and $(s,f) \in \mathcal{A}$. The set of constraints \eqref{mc_2}--\eqref{mc_5} can be reduced to
    \begin{subequations}
    \begin{align}
    & \zeta_{f,s,i} \ge 0, \quad  s\in \mathcal{F}, (s,f) \in \mathcal{A}\\
    & \zeta_{f,s,i} \ge  y_{s,i-1} -1, \quad  s\in \mathcal{F}, (s,f) \in \mathcal{A} \\
    & \zeta_{f,s,i} \le 0, \quad  s\in \mathcal{F}, (s,f) \in \mathcal{A} \\
    & \zeta_{f,s,i} \le y_{s,i-1}, \quad  s \in \mathcal{F}, (s,f) \in \mathcal{A}.
    \end{align}
    \end{subequations}
    These constraints can be further reduced to $\zeta_{f,s,i} = 0$ and $0 \le y_{s,i-1} \le 1$.
    \item When $\lambda_{f,i} = 1$, $\zeta_{f,s,i} = \lambda_{f,i} y_{s,i-1}$ gives $\zeta_{f,s,i} = y_{s,i-1}$ and $y_{s,i-1}$ are continuous between zero and one for $s \in \mathcal{F}$ and $(s,f) \in \mathcal{A}$. The set of constraints \eqref{mc_2}--\eqref{mc_5} can be reduced to
    \begin{subequations}
    \begin{align}
    & \zeta_{f,s,i} \ge 0, \quad  s\in \mathcal{F}, (s,f) \in \mathcal{A}\\
    & \zeta_{f,s,i} \ge y_{s,i-1}, \quad  s\in \mathcal{F}, (s,f) \in \mathcal{A}\\
    & \zeta_{f,s,i} \le 1, \quad  s\in \mathcal{F}, (s,f) \in \mathcal{A}\\
    & \zeta_{f,s,i} \le y_{s,i-1}, \quad  s \in \mathcal{F}, (s,f) \in \mathcal{A}.
    \end{align}    
    \end{subequations}
    These constraints can be further reduced to $\zeta_{f,s,i} = y_{s,i-1}$ and $0 \le y_{s,i-1} \le 1$.
\end{itemize}

As a result, we see that for $\lambda_{f,i}\in\{0,1\}$, the bilinear constraint $\zeta_{f,s,i} = \lambda_{f,i} y_{s,i-1}$ and the McCormick envelopes are equivalent. This completes the proof.
\endproof

\begin{corollary} 
The nonlinear follower problem \eqref{nonlinear_follower} can be equivalently written as the following mixed-binary linear program:
\begin{align*}
    \max_{\boldsymbol{y},\boldsymbol{\lambda},\boldsymbol{v}} \quad & \sum_{f \in \mathcal{F}}\sum_{i=1}^{N_p}W_f y_{f,i}	\\
    \text{subject to} \quad 
    & \eqref{inner_constr_capacity_first_direct}-\eqref{inner_constr_capacity_multistage_direct}, \eqref{nonlinear_follower:cons1}-\eqref{nonlinear_follower:cons2}, \eqref{nonlinear_follower:cons4}-\eqref{nonlinear_follower:cons6}, \eqref{mc_1}-\eqref{mc_5}\notag \\
    & \lambda_{f,i} \in \{0,1\}, \quad f \in \mathcal{F}, i \in \{1,\ldots,N_p\}.
\end{align*}
\end{corollary}


Note that the binary restrictions on $\lambda_{f,i}$ can be dealt with more efficiently in a branch-and-bound scheme with a single branching for each $\lambda_{f,i}$, whereas multiple branching may be required on the feasible region in the formulation with the bilinear constraint for each $\lambda_{f,i}$ (e.g., in a spatial branch-and-bound scheme). We will perform 
numerical experiments with both formulations for a comparison in an empirical study in Section~\ref{sec:numerical}. 
Additionally, we have observed that there exist networks where the solutions of the formulation with the McCormick envelopes and continuous $\lambda_{f,i}$ differ from those of the formulation with the bilinear constraint. As a result, the McCormick envelopes with continuous $\lambda_{f,i}$ are not exact in general. 

We conclude this section by denoting the two equivalent problem formulations. We denote the bilevel problems resulting from Proposition 1 and Corollary 1 by $(P_{nl})$ and $(P_{mc})$, respectively. 

\section{Decomposition approach} \label{sec:bender}
In this section we present a Benders-type decomposition algorithm to solve formulations  of the cascading failure problem. First, for  ease of exposition, we rewrite the problem, either $(P_{nl})$ or $(P_{mc})$, in a compact form, where, in addition to $\boldsymbol{x}$ and $\boldsymbol{y}$, we write $\boldsymbol{\lambda}$, $\boldsymbol{v}$ as the vector forms of the corresponding variables, respectively, as
\begin{subequations}
\begin{align}
    \min_{\boldsymbol{x} \in \mathcal{X}} \; \max_{\boldsymbol{y}, \boldsymbol{\lambda}, \boldsymbol{v}} \quad & W^\top \boldsymbol{y} \label{bender_compact_obj}\\
    \text{subject to} \quad & g(\boldsymbol{x}, \boldsymbol{y}, \boldsymbol{\lambda}, \boldsymbol{v}) \le 0, \label{bender_compact_constr}
\end{align}
\end{subequations}
where $\mathcal{X}$ is the feasible region for $\boldsymbol{x}$ and is the same as defined before. $g(\boldsymbol{x}, \boldsymbol{y}, \boldsymbol{\lambda}, \boldsymbol{v}) \le 0$ denote the constraints of the follower problem. We first present the master problem (MP),
\begin{subequations}
\begin{align} 
   \text{(MP)} \;\;\; \min\; & z \label{bender:master_obj}\\
    \text{s.t.} \; & z \ge Q(\boldsymbol{x}^k) - \sum_{f \in \mathcal{F}: x^k_f = 0} \alpha_f (x_f - x^k_f), \quad  k \in K\label{bender:master_optimality_cut}\\
    & \boldsymbol{x} \in \mathcal{X} \label{bender:master_x_domain},
\end{align}
\end{subequations}
where $K$ is the index set for the previously explored leader solution $\boldsymbol{x}^k$ and $Q(\boldsymbol{x}^k)$ is the objective value of the follower problem for $\boldsymbol{x}^k$. 
Constraint~\eqref{bender:master_optimality_cut} represents the optimality cuts and provides lower bounds on the objective value $z$. Note that for any $\boldsymbol{x}$, there exist feasible solutions for $\boldsymbol{y}$, $\boldsymbol{\lambda}$, and $\boldsymbol{v}$ by setting all of them to zero. The coefficients $\alpha_f$ for $f \in \mathcal{F}$ compute the reduction in the service level of the network if asset $f$ is disabled, that is, $x_f$ changes from $0$ to $1$. We compute these coefficients for $f \in \mathcal{F}$ in Algorithm \ref{algo:capacity_loss}. We use $[N_p]$ to denote the set $\{1,\ldots,N_p\}$ and $\overline{P}_{f,s,i}$ to represent upper bounds on the arc weights $P_{f,s,i}$ for $f, \; s \in \mathcal{F}$ and $i \in [N_p]$. The calculations of the upper bounds on the arc weights $\overline{P}_{f,s,i}$ depend on the description of the set $\mathcal{P}_{s,i}$. Since Algorithm \ref{algo:capacity_loss}  relies only on the network structure and number of stages of failure cascade, it  needs to be executed only once before the first iteration.

\begin{algorithm}[ht]
\caption{Computation of $\alpha_f$} \label{algo:capacity_loss}
\textbf{Input:} $W_f$ for $f \in \mathcal{F}$, $\mathcal{A}$, $\overline{P}_{f,s,i}$ for $f,\;s \in \mathcal{F}$ such that $(f,s) \in \mathcal{A}$ and $i \in \{1,\ldots, N_p\}$\\
\textbf{Output:} $\alpha_f$ for $f \in \mathcal{F}$\\
Let $\alpha_{f,i}$ be a table for $f \in \mathcal{F}$ and $i \in [N_p]$\\
\For{$f \in \mathcal{F}$}{
$\alpha_{f,1} = W_f + \sum_{s \in \mathcal{F}: (f,s) \in \mathcal{A}} W_s\overline{P}_{f,s,1}$\\
\For{$i \in \{2,\ldots,N_p\}$}{
$\alpha_{f,i} = W_f + \sum_{s \in \mathcal{F}: (f,s) \in \mathcal{A}} \alpha_{s,i-1} \overline{P}_{f,s,i}$
}
}
\For{$f \in {\cal F}$}{
$\alpha_f = \sum_{i \in [N_p]} \alpha_{f,i}$
}
\Return{$\alpha_f$ for $f \in \mathcal{F}$}
\end{algorithm}

Next, from the construction of the master problem, the subproblem is exactly the follower problem, $Q(\boldsymbol{x}^*)$, to solve for the maximum remaining service level of the network given the leader solution $\boldsymbol{x}^*$. With the same notation from the compact form~\eqref{bender_compact_obj}--\eqref{bender_compact_constr}, the subproblem (SP) is given by
\begin{subequations}
\begin{align}
    (\text{SP}) \;\;\; Q(\boldsymbol{x}^*) = & \max_{\boldsymbol{y}, \boldsymbol{\lambda}, \boldsymbol{v}} W^\top \boldsymbol{y} \label{bender:sub_obj}\\
    & \text{s.t.} \; g(\boldsymbol{x}^*, \boldsymbol{y}, \boldsymbol{\lambda}, \boldsymbol{v}) \le 0. \label{bender:sub_constr}
\end{align}
\end{subequations}

We are now ready to present the overall procedure in Algorithm \ref{algo:bender} with a target accuracy $\varepsilon$ in percentage. We denote an upper bound for $z$ by $\overline{z}$ and initialize it to infinity, while we denote a lower bound for $z$ by $\underline{z}$ and initialize it to zero as $\boldsymbol{y} \ge 0$. In practice, we can use any feasible solution of the master problem  to construct a subproblem and obtain a cut. However, only the objective values associated with optimal solutions are valid lower bounds to update $\underline{z}$ in Step \ref{update_lower_bound} of Algorithm \ref{algo:bender}. 

\begin{algorithm}[ht]
\caption{Decomposition algorithm} \label{algo:bender}	
\textbf{Input:} Target gap $\varepsilon \ge 0$\\
Initialize upper bound $\overline{z} = \infty$, lower bound $\underline{z} = 0$, and iteration count $k = 0$\\
Compute the cut coefficients $\alpha_f$\\
\While{$\overline{z} \ge (1 + \varepsilon) \underline{z}$}{
$k \leftarrow k + 1$\\
Solve the master problem (MP) to obtain the solution $\boldsymbol{x}^k$ and $z^k$\\
\lIf{$z^k \ge \underline{z}$}{update $\underline{z} = z^k$} \label{update_lower_bound} \label{underline_z_update}
Fix $\boldsymbol{x} = \boldsymbol{x}^k$ and solve the subproblem (SP)\\
\lIf{$Q(\boldsymbol{x}^k) \le \overline{z}$}{update $\overline{z} = Q(\boldsymbol{x}^k)$}
Add cut~\eqref{bender:master_optimality_cut} to (MP)\\
}
\Return{$\varepsilon$-optimal solution $x^k$}
\end{algorithm}

Next, we present a modification of Algorithm~\ref{algo:bender}.
We consider the strengthening of cut~\eqref{bender:master_optimality_cut} in the form of 
\begin{equation}
	z \ge Q(\boldsymbol{x}^k) - \sum_{f \in \mathcal{F}: x^k_f = 0} \alpha_f (x_f - x^k_f) + \sum_{f \in \mathcal{F}: x^k_f = 1}\beta_f(\boldsymbol{x}^k)(x_f - x^k_f) \label{bender:master_optimality_cut_strengthen}.
\end{equation}
The coefficients $\beta_f(\boldsymbol{x}^k)$ compute the release of service level of the network if an asset $f$ that is disabled in $\boldsymbol{x}^k$ is instead operational. Unlike the coefficients of $\alpha_f$, $\beta_f(\boldsymbol{x}^k)$ depends on a particular leader solution $\boldsymbol{x}^k$, and we can utilize auxiliary problems to compute them. 
Consider a particular $\boldsymbol{x}^k$ and $f$ such that $x^k_f = 1$; we use $\boldsymbol{\tilde{x}}$ to denote the new vector obtained from changing $x^k_f$ to $0$ while keeping the rest of components unchanged. We solve the subproblem $Q(\boldsymbol{\tilde{x}})$ and set 
\begin{equation}
    \beta_f(\boldsymbol{x}^k) = Q(\boldsymbol{\tilde{x}}) - Q(\boldsymbol{x}^k).
\end{equation}
We repeat the process for all assets that are disabled in $\boldsymbol{x}^k$ to obtain a strengthened cut~\eqref{bender:master_optimality_cut_strengthen}. Depending on the value of $N_c$, this may require solving a few auxiliary problems. Note that we can add cut~\eqref{bender:master_optimality_cut} for each $\tilde{x}$ as well.

We now discuss the convergence of Algorithm~\ref{algo:bender}, which is evident after first developing Lemma~\ref{lemma:convergence}. Lemma~\ref{lemma:convergence} shows that if, in any iteration, the optimal solution is the same as any of the previous optimal solutions, then Algorithm~\ref{algo:bender} terminates.
Formally, we have the following.
\begin{lemma} \label{lemma:convergence}
If, for any iteration $n+1$, the optimal solution from solving (MP) $\boldsymbol{x}^{n+1} = \boldsymbol{x}^j$ for some $j \in \{1,\ldots,n\}$, then $\overline{z} = \underline{z}$.  
\end{lemma}
\proof{Proof of Lemma \ref{lemma:convergence}.} 
In iteration $n+1$, we have that the optimal objective value 
\begin{equation}
    z^{n+1} \ge \underline{z} = \max\{z^1,\ldots,z^n\}. 
\end{equation}
The first inequality is due to the fact that (MP) in iteration $n+1$ has more constraints than any previous iteration does. Additionally, by the definition of $\overline{z}$, we have that
\begin{equation}
    \overline{z} = \min\{Q(\boldsymbol{x}^1), \ldots, Q(\boldsymbol{x}^n)\}. \label{overline_z_def}
\end{equation}
Suppose that $\boldsymbol{x}^{n+1} = \boldsymbol{x}^j$ for some $j \in \{1,\ldots,n\}$. Then the optimality cut for $\boldsymbol{x}^j$, either as \eqref{bender:master_optimality_cut} or as \eqref{bender:master_optimality_cut_strengthen}, reduces to
\begin{equation}
    z^{n+1} \ge Q(\boldsymbol{x}^j).
\end{equation}
Since $z^{n+1}$ and $\overline{z}$ are a lower bound and an upper bound to the original problem, respectively, we have that $\overline{z} \ge z^{n+1}$. Additionally, \eqref{overline_z_def} implies $Q(\boldsymbol{x}^j) \ge \overline{z}$. As a result, we have the following inequality,
\begin{equation}
   Q(\boldsymbol{x}^j) \ge \overline{z} \ge z^{n+1} \ge Q(\boldsymbol{x}^j).
\end{equation}
Thus, we have $\overline{z} = Q(\boldsymbol{x}^j)$. Next, the update step for $\underline{z}$ (see line \ref{underline_z_update} of Algorithm \ref{algo:bender}) implies
\begin{equation}
    \underline{z} = \max\{\underline{z}, Q(\boldsymbol{x}^j)\} \ge Q(\boldsymbol{x}^j) = \overline{z}.
\end{equation}
\endproof
With Lemma \ref{lemma:convergence}, we are ready to state the convergence of the algorithm. 
\begin{theorem} \label{thm:finite_converge}
Algorithm \ref{algo:bender} terminates after a finite number of iterations. In particular, the algorithm terminates after exploring all $\sum_{n \le N_c} \binom{N}{n}$ feasible leader solutions $\boldsymbol{x}$ in the worst case. 
\end{theorem}
\proof{Proof of Theorem \ref{thm:finite_converge}.}
Recall the budget constraint~\eqref{constr_budget} for the leader $e^\top \boldsymbol{x} \le N_c$, which implies that there are in total $\sum_{n \in N_c} \binom{N}{n}$ feasible leader solutions $\boldsymbol{x}$. As a result of Lemma \ref{lemma:convergence}, the algorithm terminates when a previous optimal solution is obtained as the optimal solution from solving (MP) in a later iteration. Consequently, Algorithm \ref{algo:bender} terminates after exploring all the feasible solutions in the worst case because the optimal solution in iteration $\sum_{n \le N_c} \binom{N}{n} + 1$ must have been explored in a previous iteration.
\endproof

\section{Empirical study: resiliency in Puerto Rico infrastructure} \label{sec:numerical}
In this section we perform numerical experiments to validate the proposed algorithms on infrastructure networks derived with anonymized real-world data. The data that was collected as part of the DIS team's support for recovery efforts in Puerto Rico consisted of both geographic information system  data and various contextual information on the assets and systems provided by the infrastructure owners and operators in Puerto Rico. This included data on the characteristics, locations, relative proximities, service areas, and customer bases of the infrastructure. The DIS team used this data to assist decision-makers in prioritizing recovery and hazard mitigation investments based on which assets and system components are the most widely required in order to enable dependent downstream operations. Because of  sensitivities of this data, a subset was produced that did not include any actual identifications of assets or system components in terms of name, owner, or location. This allowed for numerical experiments on anonymized data that nevertheless presented the characteristics, relative proximities, and service areas of the real-world data that are helpful in understanding how these infrastructures connect to and depend on other infrastructure. 
\subsection{Input data and the networks}
We use the anonymized real-world data to generate infrastructure networks to validate the proposed algorithm. The networks are generated based on the frequency or distribution of the asset classes from Table~\ref{table:asset_depedency}. The importance weight, $W$, for each asset class is assumed to be inversely proportional to its frequency. We also assume that assets from the same asset class share the same importance weights. We present the frequency and importance weights in Table \ref{table:asset_info}.
\begin{table}[ht]
\caption{Asset class frequency and weights} \label{table:asset_info}
	\centering
    \small
	\begin{tabular}{|c|c|c|}
	\hline
	Asset class  & Frequency & $W$ \\\hline\hline
	Cell tower & 0.0438 & 23 \\\hline
	EMS & 0.0501 & 20 \\\hline
	Thermal generation & 0.00887 & 113 \\\hline
	Renewable generation & 0.0162 & 62 \\\hline
	Substation & 0.178 & 6 \\\hline
	Health care & 0.0391 & 26 \\\hline
	Pharmacies & 0.527 & 2 \\\hline
	Transport & 0.0318 & 31 \\\hline
	Water & 0.105 & 10 \\\hline
	\end{tabular}
\end{table}
After generating the assets according to the frequencies, we describe how to generate the set $\mathcal{P}_f$ for $f \in \mathcal{F}$. For each $f \in \mathcal{F}$, we consider any asset $s$ such that $(s,f) \in \mathcal{A}$. If we denoted the Euclidean distance between these two assets by $d_{s,f}$, we use a Hoffman model (inverse distance weighting in~\cite{Shepard1968}) to compute a base value for the dependencies, denoted by $\tilde{P}_{s,f}$, for all $i \in \{1,\ldots,N_p\}$, which are given by
\begin{equation}
\tilde{P}_{s,f} = \frac{1}{d^2_{s,f}},\; \forall (s,f) \in \mathcal{A}.	
\end{equation}
We set a threshold value for the distance $d_{s,f}$ below which we set $\tilde{P}_{s,f}$ to be zero since those assets are too far away from each other to have a dependency. For asset $s$ such that $(s,f) \notin \mathcal{A}$, we also set the base value $\tilde{P}_{s,f}$ to zero. We then consider a set $\mathcal{P}_f$ in our numerical experiments as follows: 
\begin{equation}
\mathcal{P}_f = \left\{P_{s,f,i} \ge 0: \; \sum_{s \in \mathcal{F}:(s,f) \in \mathcal{A}} P_{s,f,i} = 1, \; |P_{s,f,i} - \tilde{P}_{s,f}| \le \delta \tilde{P}_{s,f}, \; \forall s \in \mathcal{F}\right\},
\end{equation}
with $\delta = 0.1$. Consequently, we can obtain the upper bounds, $\overline{P}_{s,f,i}$, for $P_{s,f,i}$, that are used in Algorithm \ref{algo:capacity_loss} as $(1 + \delta)\tilde{P}_{s,f}$. 

For our experiments, we generate two networks that vary in the number of assets. 
Table \ref{table:network_info} shows a breakdown on the size of the networks. 
\begin{table}[ht]
	\caption{Network information} \label{table:network_info}
	\centering
    \small
	\begin{tabular}{|c|c|c|}
	\hline
		Network name & \texttt{a128} & \texttt{a210} \\\hline\hline
		No. of assets & 128 & 210  \\\hline
		No. of arcs in $\mathcal{A}$ & 1044 & 2033 \\\hline
	\end{tabular}
\end{table}
We  point out that many larger infrastructure networks can be mapped into smaller networks of sizes that are comparable to our infrastructure networks by grouping infrastructure assets by their geographical proximity. This is valid because often natural disaster events can disable multiple assets located in an area rather than a single asset. The dependency graph can be constructed by considering the functionalities of the assets in the groups. As a result, we are interested in identifying the most important asset groups in the networks.

\subsection{Computational settings}
We run the experiments on a computer with an Intel i9 CPU with 64 GB RAM. The computer runs the Ubuntu 20.04 LTS operating system. The algorithm is implemented in Julia 1.6.1 (\cite{Julia-2017}) and JuMP. We use Gurobi 9.5.1 (\cite{gurobi}) to solve both the master problem  and subproblem. For each instance, we set an 1-hour time limit. We set the target gap to be $1\%$ while using $10$ feasible solutions including the optimal solution from the master problem  to generate optimality cuts. 

\subsection{Computational performances}
In this section we compare the performances between the formulations $(P_{nl})$ and $(P_{mc})$ as well as across the variants of the algorithm. We denote Algorithm \ref{algo:bender} with optimality cuts~\eqref{bender:master_optimality_cut} by \texttt{Benders} and with strengthened optimality cuts~\eqref{bender:master_optimality_cut_strengthen} by \texttt{Benders-S}, respectively. Since strengthening requires solving auxiliary problems, it is not practical to perform strengthening for all feasible solutions, and we  strengthen only the optimality cut for the optimal solution from the master problem.

\subsubsection{Comparing the formulations $(P_{nl})$ and $(P_{mc})$}
We first compare the performances between the formulations $(P_{nl})$ and $(P_{mc})$. For $(P_{nl})$, we take advantage of the advancements for quadratic programs in Gurobi 9.5.1 and solve it directly with the solver. The detailed results are reported in Table \ref{table:nc_mc}. The column ``Iter." reports the number of iterations. 

\begin{table}[ht]
\caption{Comparison of formulations $(P_{nl})$ and $(P_{mc})$}	 \label{table:nc_mc}
\centering
\small
\begin{tabular}{|c|c|c|c|c|c|c|c|}
	\hline
	Instance & Formulation \& method & $N_p$ & $N_c$ & Time (s.) & Iter. & Opt. obj & Disabled assets \\\hline
	\multirow{4}{*}{\texttt{a128}} & $(P_{nl})$ \& \texttt{Benders} & 3 & 5 & 176.87 & 24 & 3754.69 & [1, 17, 33, 67, 97]\\\cline{2-8}
    & $(P_{nl})$ \& \texttt{Benders-S} & 3 & 5 & 168.14 & 15 & 3754.69 & [1, 17, 33, 67, 97]\\\cline{2-8}
	& $(P_{mc})$ \& \texttt{Benders} & 3 & 5 & 57.48 & 24 & 3754.69 & [1, 17, 33, 67, 97]\\\cline{2-8}
    & $(P_{mc})$ \& \texttt{Benders-S} & 3 & 5 & 50.16 & 15 & 3754.69 & [1, 17, 33, 67, 97]\\\hline
	\multirow{4}{*}{\texttt{a210}} & $(P_{nl})$ \& \texttt{Benders} & 3 & 5 & 520.95 & 30 & 6052.88 & [1, 43, 85, 106, 148]\\\cline{2-8}
    & $(P_{nl})$ \& \texttt{Benders-S} & 3 & 5 & 425.75 & 17 & 6052.88 & [1, 43, 85, 106, 148]\\\cline{2-8}
    & $(P_{mc})$ \& \texttt{Benders} & 3 & 5 & 166.89 & 29 & 6052.88 & [1, 43, 85, 106, 148]\\\cline{2-8}
	& $(P_{mc})$ \& \texttt{Benders-S} & 3 & 5 & 138.67 & 17 & 6052.88 & [1, 43, 85, 106, 148]\\\hline
\end{tabular}
\end{table}

In this set of experiments we choose the number of stages of failure cascade ($N_p$) and number of disabled assets ($N_c$) to be such that both the $(P_{nl})$ and $(P_{mc})$ can be solved within the time limit. We observe that the $(P_{nl})$ and $(P_{mc})$ yield the same optimal objective values and solutions (i.e., disabled assets). This set of results empirically confirms  the exactness proof in Theorem~\ref{theorem:mccormick_equiv}. Nonetheless, the computation times with the $(P_{mc})$ are significantly shorter. The average improvement in computation time from $(P_{mc})$ is about $68.17\%$, while the number of iterations needed to achieve the target gap does not differ significantly for both formulations. Additionally, \texttt{Benders-S} demonstrates further improvements from \texttt{Benders}, benefiting from the strengthened optimality cuts~\eqref{bender:master_optimality_cut_strengthen}.

\subsubsection{Comparing the variants \texttt{Benders} and \texttt{Benders-S}}
Next, we compare the performances between the variants \texttt{Benders} and \texttt{Benders-S}. In this set of experiments we use the $(P_{mc})$. We first present the plots for the progress of upper bound $\overline{z}$ and the lower bound $\underline{z}$ as well as the gaps as computed by 
\begin{equation}
    \text{Gap} = \frac{\overline{z} - \underline{z}}{\underline{z}} \times 100\%.
\end{equation}

\begin{figure*}[t!]
    \centering
    \begin{minipage}[t]{0.5\textwidth}
        \centering
        \includegraphics[width = \linewidth]{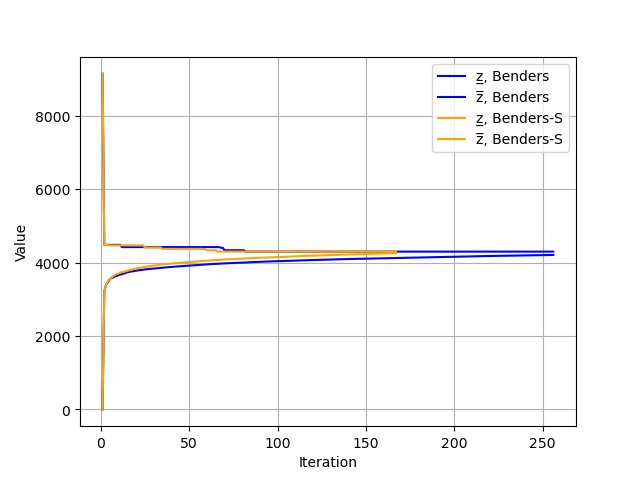}
        \caption{Full progress of $\overline{z}$ and $\underline{z}$ for the network \texttt{a210} with $N_p = 3$ and $N_c = 8$}\label{fig:full_progress}
    \end{minipage}%
    ~ 
    \begin{minipage}[t]{0.5\textwidth}
        \centering
        \includegraphics[width = \linewidth]{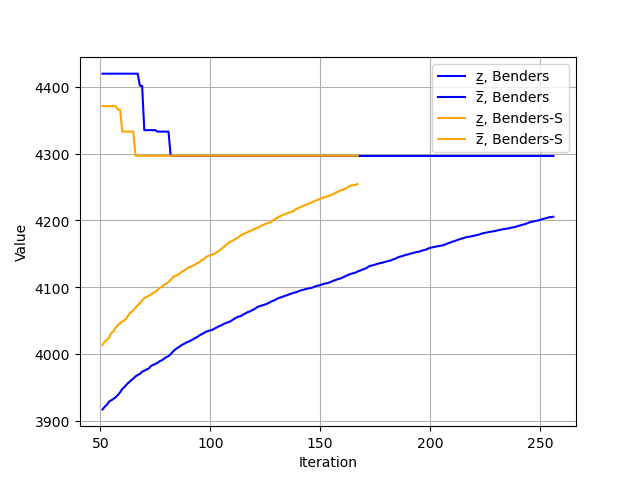}
        \caption{Progress of $\overline{z}$ and $\underline{z}$ from iteration $50$ onwards for the network \texttt{a210} with $N_p = 3$ and $N_c = 8$}\label{fig:zoom_progress}
    \end{minipage}
\end{figure*}

\begin{figure*}[t!]
    \centering
    \begin{minipage}[t]{0.5\textwidth}
        \centering
        \includegraphics[width = \linewidth]{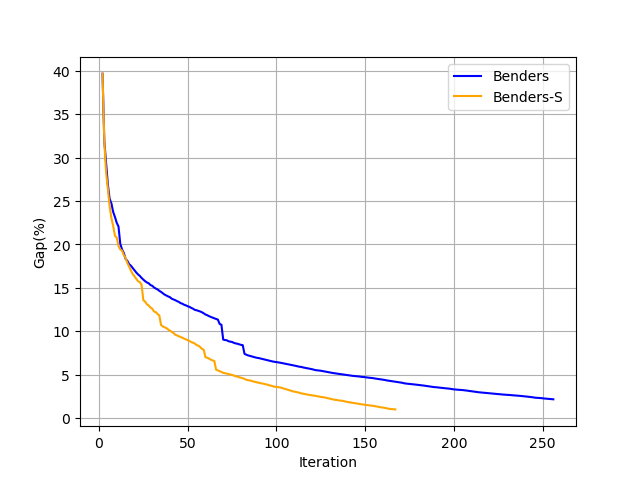}
        \caption{Progress of gap for the network \texttt{a210} with $N_p = 3$ and $N_c = 8$}\label{fig:gaps}
    \end{minipage}%
    ~ 
    \begin{minipage}[t]{0.5\textwidth}
        \centering
        \includegraphics[width = \linewidth]{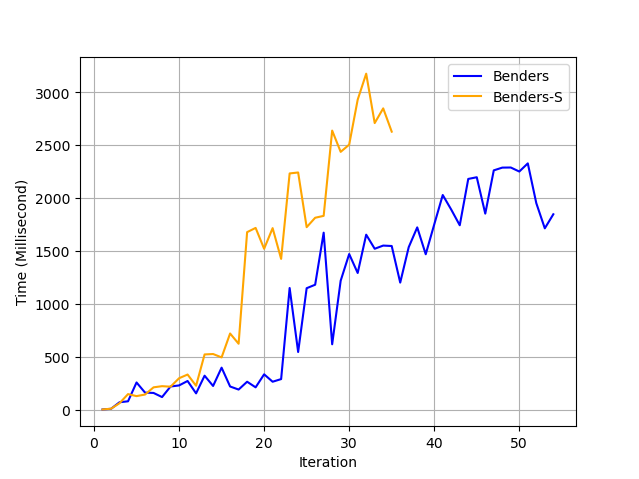}
        \caption{Time taken to solve the master problem (MP) for the network \texttt{a210} with $N_p = 3$ and $N_c = 8$}\label{fig:MP_time}
    \end{minipage}
\end{figure*}

Figures~\ref{fig:full_progress}--\ref{fig:gaps} show the optimality gaps of various solvers for the network \texttt{a210} with 3 stages of failure cascade ($N_p = 3$) and 8 disabled assets ($N_c = 8$). We observe that the gap between $\overline{z}$ and $\underline{z}$ reduces  quickly in the first $30$ iterations for both variants. After approximately 80 iterations, the upper bound $\overline{z}$ remains unchanged while computations are carried out to close the gaps through improving the lower bound $\underline{z}$. This observation suggests that the optimal solution and objective values are likely to be found early in the execution of the algorithm, while proving optimality is a much more challenging task. In addition, the number of iterations needed to reach the target gap for \texttt{Benders-S} is  significantly smaller than that for \texttt{Benders}. We point out that \texttt{Benders} terminates with a $2.16\%$ gap, whereas \texttt{Benders-S} achieves the target gap in the given time limit. This result shows the effectiveness of the strengthened optimality cuts despite the fact that, in each iteration, additional auxiliary problems have to be solved in the \texttt{Benders-S} variant.

We also track the time taken to solve the master problem ~\eqref{bender:master_obj}-~\eqref{bender:master_x_domain} in each iteration for both \texttt{Benders} and \texttt{Benders-S}. Figure~\ref{fig:MP_time} shows a plot for the network \texttt{a128} with $2$ stages of failure cascade $(N_p = 2)$ and $9$ disabled assets $(N_c = 9)$. We observe a general trend of longer time taken as the number of iterations increases for both \texttt{Benders} and \texttt{Benders-S}. This is expected as optimality cuts are added in each iteration to the master problem, increasing the size of the problem. The rate of increment, however, is larger for \texttt{Benders-S}  than it is for \texttt{Benders}. This is due to the additions of cuts~\eqref{bender:master_optimality_cut} for each $\boldsymbol{\tilde{x}}$ in the strengthening process, as discussed in Section~\ref{sec:bender}.


\subsubsection{Scalability of the algorithm}
The next set of results probe the scalability of our algorithm. 
In particular, we vary the number of stages of failure cascade ($N_p$) and number of disabled assets ($N_c$) for both networks. Similar to the previous set of experiments, we use $(P_{mc})$. In Tables \ref{table:a128} and \ref{table:a210} we report the detailed computational results. From both tables we see that with a smaller number of stages of failure cascade ($N_p = 1$ or $2$), both variants of the algorithms solve the problem to optimality in a relatively short amount of time. As we increase the number of stages of failure cascade, the problems become much more challenging to solve. When we consider 4 stages of failure cascade ($N_p = 4$), in few cases we reach the time limit with non-zero gaps. Nonetheless, we observe similar progression of the values of $\overline{z}$ and $\underline{z}$ as well as gaps to what are shown in Figures \ref{fig:full_progress}-- \ref{fig:zoom_progress}. The algorithm is able to find feasible solutions in these cases, but more computations are needed to improve the lower bounds to close the gap and prove that the feasible solutions are optimal. Even with the non-zero gaps, these feasible solutions have their practical uses. In addition, we see that increasing the number of stages of failure cascade ($N_p$) has a more significant impact than increasing the number of disabled assets ($N_c$) on the computation. This observation aligns with our expectation. As we increase the number of stages of failure cascade ($N_p$), the follower problem (SP) becomes larger with more variables for the expected fractional service level; in contrast,  as we increase the number of disabled assets ($N_c$), we do not change the size of the follower problem (SP). When the problems are relatively easy with small number of stages of failure cascade ($N_p$) or disabled assets ($N_c$), we see that \texttt{Benders} outperforms \texttt{Benders-S}, because \texttt{Benders-S} requires solving additional auxiliary problems in the process. The effectiveness and benefits of the strengthened optimality cuts become obvious when the problems are challenging to solve with larger numbers of stages of failure cascade ($N_p$) or disabled assets ($N_c$). In a few instances  \texttt{Benders} terminates with non-zero gaps while \texttt{Benders-S} is able to close the gaps despite longer time taken per  iteration with \texttt{Benders-S}.

\begin{table}[ht!]
\caption{Computational results for the network \texttt{a128}}\label{table:a128}
\centering
\small
\begin{tabular}{|c|c|c|c|c|c|c|c|}
\hline
	\multirow{2}{*}{$N_p$} & \multirow{2}{*}{$N_c$} & \multicolumn{3}{c|}{\texttt{Benders}} & \multicolumn{3}{c|}{\texttt{Benders-S}} \\\cline{3-8}
	& & Time (s.) & Iter. & Gap(\%) & Time (s.) & Iter. & Gap(\%) \\\hline\hline
	\multirow{10}{*}{1} & 2 & 1.48 & 2 & 0.57 & 1.48 & 2 & 0.57\\\cline{2-8}  
	& 3 & 1.97 & 3 & 0.0 & 1.95 & 3 & 0.0 \\\cline{2-8}
	& 4 & 1.89 & 3 & 0.0 & 2.06 & 3 & 0.0 \\\cline{2-8}
    & 5 & 1.98 & 3 & 0.26 & 2.10 & 3 & 0.0 \\\cline{2-8}
    & 6 & 2.28 & 4 & 0.79 & 2.71 & 4 & 0.11 \\\cline{2-8}
    & 7 & 2.93 & 6 & 0.73 & 3.54 & 5 & 0.68 \\\cline{2-8}
    & 8 & 5.80 & 12 & 0.82 & 5.79 & 8 & 0.99 \\\cline{2-8}
    & 9 & 7.45 & 16 & 0.52 & 8.69 & 11 & 0.03 \\\cline{2-8}
	& 10 & 9.52 & 19 & 0.73 & 9.04 & 11 & 0.75 \\\hline\hline
	\multirow{10}{*}{2} & 2 & 3.36 & 2 & 0.99 & 3.52 & 2 & 0.99\\\cline{2-8}
	& 3 & 4.63 & 3 & 0.0 & 4.97 & 3 & 0.0\\\cline{2-8}
	& 4 & 5.88 & 4 & 0.71 & 7.03 & 4 & 0.0\\\cline{2-8}
    & 5 & 7.23 & 5 & 0.9 & 9.47 & 5 & 0.48\\\cline{2-8}
    & 6 & 11.17 & 8 & 0.83 & 12.54 & 6 & 0.65 \\\cline{2-8}
    & 7 & 23.89 & 16 & 0.86 & 23.58 & 11 & 0.68 \\\cline{2-8}
    & 8 & 44.58 & 27 & 0.86 & 42.87 & 18 & 0.77 \\\cline{2-8}
    & 9 & 123.89 & 54 & 0.88 & 121.96 & 35 & 0.92 \\\cline{2-8}
	& 10 & 356.62 & 96 & 0.95 & 340.83 & 60 & 0.91\\\hline\hline
	\multirow{10}{*}{3} & 2 & 6.95 & 3 & 0.0 & 7.94 & 3 & 0.0\\\cline{2-8}
	& 3 & 9.37 & 4 & 0.0 & 8.47 & 3 & 0.0\\\cline{2-8}
	& 4 & 20.22 & 9 & 0.77 & 20.03 & 7 & 0.42\\\cline{2-8}
    & 5 & 57.48 & 24 & 0.86 & 50.16 & 15 & 0.71\\\cline{2-8}
    & 6 & 107.68 & 38 & 0.83 & 96.23 & 25 & 0.95 \\\cline{2-8}
    & 7 & 252.17 & 67 & 0.95 & 229.51 & 43 & 0.86 \\\cline{2-8}
    & 8 & 664.58 & 121 & 0.99 & 637.83 & 77 & 0.86 \\\cline{2-8}
    & 9 & 2528.53 & 243 & 0.98 & 1971.10 & 150 & 0.97 \\\cline{2-8}
	& 10 & 3600 & 237 & 9.40 & 3600 & 190 & 5.89\\\hline\hline
	\multirow{10}{*}{4} & 2 & 11.86 & 3 & 0.0 & 13.52 & 3 & 0.0\\\cline{2-8}
	& 3 & 22.65 & 6 & 0.0 & 22.40 & 5 & 0.57\\\cline{2-8}
	& 4 & 72.57 & 19 & 0.56 & 72.51 & 14 & 0.34\\\cline{2-8}
    & 5 & 178.65 & 42 & 0.86 & 171.92 & 29 & 0.95\\\cline{2-8}
    & 6 & 486.85 & 87 & 0.99 & 463.02 & 60 & 0.88 \\\cline{2-8}
    & 7 & 1219.41 & 144 & 0.96 & 1015.32 & 96 & 0.92 \\\cline{2-8}
    & 8 & 3600 & 282 & 2.04 & 3126.16 & 196 & 0.96 \\\cline{2-8}
    & 9 & 3600 & 256 & 13.47 & 3600 & 189 & 11.79 \\\cline{2-8}
	& 10 & 3600 & 204 & 55.03 & 3600 & 145 & 53.85\\\hline	
\end{tabular}
\end{table}

\begin{table}[ht!]
\caption{Computational results for the network \texttt{a210}}\label{table:a210}
\centering
\small
\begin{tabular}{|c|c|c|c|c|c|c|c|}
\hline
	\multirow{2}{*}{$N_p$} & \multirow{2}{*}{$N_c$} & \multicolumn{3}{c|}{\texttt{Benders}} & \multicolumn{3}{c|}{\texttt{Benders-S}} \\\cline{3-8}
	& & Time (s.) & Iter. & Gap(\%) & Time (s.) & Iter. & Gap(\%)\\\hline\hline
	\multirow{10}{*}{1} & 2 & 4.95 & 2 & 0.23 & 4.92 & 2 & 0.23\\\cline{2-8}  
	& 3 & 4.93 & 2 & 0.79 & 5.02 & 2 & 0.79\\\cline{2-8}
	& 4 & 5.64 & 3 & 0.0 & 6.34 & 3 & 0.0 \\\cline{2-8}
    & 5 & 5.72 & 3 & 0.0 & 6.29 & 3 & 0.0 \\\cline{2-8}
    & 6 & 5.62 & 3 & 0.14 & 6.50 & 3 & 0.0 \\\cline{2-8}
    & 7 & 6.56 & 4 & 0.46 & 6.52 & 3 & 0.0 \\\cline{2-8}
    & 8 & 8.11 & 6 & 0.70 & 11.42 & 6 & 0.92 \\\cline{2-8}
    & 9 & 13.73 & 12 & 0.97 & 16.38 & 9 & 0.87 \\\cline{2-8}
	& 10 & 27.80 & 25 & 0.96 & 28.64 & 15 & 0.90 \\\hline\hline
	\multirow{10}{*}{2} & 2 & 9.84 & 2 & 0.42 & 11.37 & 2 & 0.42 \\\cline{2-8}
	& 3 & 12.67 & 3 & 0.0 & 14.52 & 3 & 0.0\\\cline{2-8}
	& 4 & 12.83 & 3 & 0.0 & 14.91 & 3 & 0.0\\\cline{2-8}
    & 5 & 12.80 & 3 & 0.75 & 15.76 & 3 & 0.46 \\\cline{2-8}
    & 6 & 18.64 & 5 & 0.48 & 20.81 & 4 & 0.11 \\\cline{2-8}
    & 7 & 21.30 & 6 & 0.95 & 26.05 & 5 & 0.96 \\\cline{2-8}
    & 8 & 94.57 & 30 & 0.94 & 56.74 & 11 & 0.96 \\\cline{2-8}
    & 9 & 302.84 & 70 & 0.99 & 280.75 & 44 & 0.93 \\\cline{2-8}
	& 10 & 733.56 & 125 & 0.97 & 615.60 & 71 & 0.97\\\hline\hline
	\multirow{10}{*}{3} & 2 & 20.24 & 3 & 0.0 & 21.59 & 3 & 0.0\\\cline{2-8}
	& 3 & 25.38 & 4 & 0.62 & 29.60 & 4 & 0.0\\\cline{2-8}
	& 4 & 53.04 & 9 & 0.95 & 53.94 & 7 & 0.67 \\\cline{2-8}
    & 5 & 166.89 & 29 & 0.92 & 138.67 & 17 & 0.85 \\\cline{2-8}
    & 6 & 379.95 & 61 & 0.98 & 352.28 & 38 & 0.91 \\\cline{2-8}
    & 7 & 965.22 & 124 & 0.98 & 760.02 & 68 & 0.96 \\\cline{2-8}
    & 8 & 3600 & 256 & 2.16 & 3389.23 & 167 & 0.99 \\\cline{2-8}
    & 9 & 3600 & 241 & 5.95 & 3600 & 165 & 5.47 \\\cline{2-8}
	& 10 & 3600 & 219 & 10.31 & 3600 & 155 & 9.65\\\hline\hline
	\multirow{10}{*}{4}& 2 & 38.82 & 4 & 0.76 & 32.05 & 3 & 0.99\\\cline{2-8}
	& 3 & 105.61 & 12 & 0.81 & 88.17 & 8 & 0.68\\\cline{2-8}
	& 4 & 452.95 & 51 & 0.93 & 285.50 & 24 & 0.89 \\\cline{2-8}
    & 5 & 1092.11 & 108 & 0.99 & 912.51 & 64 & 0.95 \\\cline{2-8}
    & 6 & 2510.36 & 188 & 0.99 & 1964.08 & 107 & 0.99 \\\cline{2-8}
    & 7 & 3600 & 229 & 4.23 & 3600 & 157 & 3.27 \\\cline{2-8}
    & 8 & 3600 & 213 & 14.24 & 3600 & 145 & 13.49 \\\cline{2-8}
    & 9 & 3600 & 196 & 22.50 & 3600 & 138 & 21.61 \\\cline{2-8}
	& 10 & 3600 & 183 & 39.22 & 3600 & 120 & 38.52\\\hline	
\end{tabular}
\end{table}

\subsection{Analyses of the network \texttt{a128}}
In this section we analyze the network \texttt{a128} for the most severe contingencies, the optimal arc weights, and demonstrate the cascading failure effects propagating through the network under certain combinations of numbers of disabled assets ($N_c$) and stages of failure cascades ($N_p$). These analyses can validate and give insights into our modeling choices. 
\subsubsection{Contingencies}
In Table~\ref{table:disabled_assets} we compare the most severe contingencies, namely, the optimal attack plan on assets, with the network \texttt{a128} under a range of values for the number of disabled assets ($N_c$) and stages of failure cascades ($N_p$). We observe that for the same number of stages of failure cascade, the most severe contingencies for smaller numbers of disabled assets ($N_c$) are not always a subset of the most severe contingencies of larger numbers of disabled assets ($N_c$). This observation aligns with what is discussed in~\cite{Morton2007}. 

\begin{table}[ht!]
\caption{Disabled assets across different $N_c$ for the network \texttt{a128}}\label{table:disabled_assets}
\centering
\small
\begin{tabular}{|c|c|c|}
\hline
	\diagbox{$N_c$}{$N_p$} & 2 & 3\\\hline
    2 & [1, 33] & [1, 33]\\\hline
    3 & [1, 30, 33] & [1, 33, 97] \\\hline 
    4 & [1, 17, 30, 33] & [1, 17, 33, 97] \\\hline
    5 & [1, 17, 30, 33, 97] & [1, 17, 33, 67, 97] \\\hline
    6 & [1, 17, 30, 33, 78, 97] & [1, 17, 30, 33, 67, 97] \\\hline
    7 & [1, 17, 30, 33, 65, 78, 97] & [1, 17, 30, 33, 67, 81, 97] \\\hline
    8 & [1, 17, 30, 33, 65, 78, 81, 97] & [1, 17, 30, 33, 67, 78, 81, 97] \\\hline
    9 & [1, 17, 30, 33, 65, 78, 81, 97, 113] & [1, 17, 30, 33, 49, 65, 78, 81, 97]
\\\hline
\end{tabular}
\end{table}

\subsubsection{Optimal arc weights}
\begin{figure*}[t!]
    \centering
    \begin{minipage}[t]{0.5\textwidth}
        \centering
        \includegraphics[width = \linewidth]{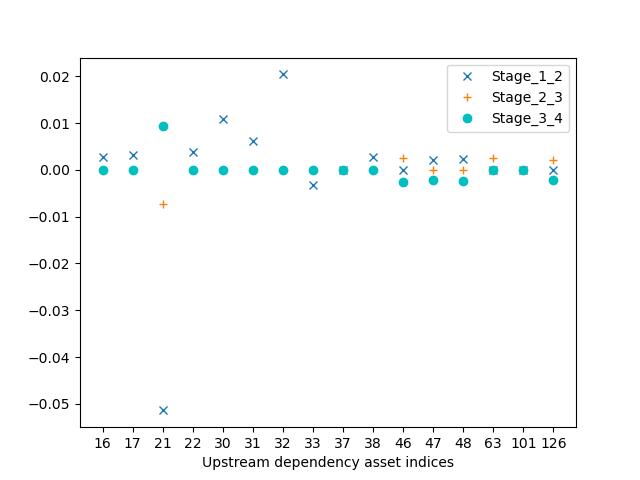}
        \caption{Changes of optimal arc weights for asset 25 ($f = 25$)}\label{fig:weights_change_1}
    \end{minipage}%
    ~ 
    \begin{minipage}[t]{0.5\textwidth}
        \centering
        \includegraphics[width = \linewidth]{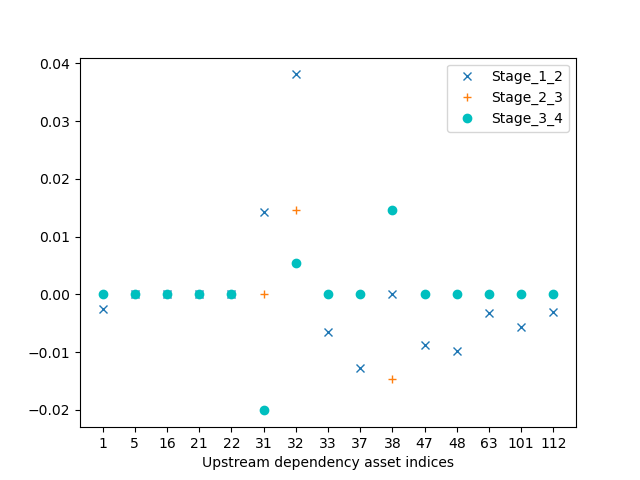}
        \caption{Changes of optimal arc weights for asset 30 ($f = 30$)}\label{fig:weights_change_2}
    \end{minipage}
\end{figure*}

Next, we compare the optimal arc weights across the stages of failure cascades' in other words, we compare how the optimal solution $P_{s,f,i}^*$ varies for $i \in \{1,\ldots,N_p\}$ for a particular $f \in \mathcal{F}$ and all $s \in \mathcal{F}$ such that $(s,f) \in \mathcal{A}$. As we discussed in Section~\ref{sec:model}, the worst-case realization $P_{s,f,i}$ can vary across the stages of failure cascade. In particular, we consider the network \texttt{a128} with $4$ stages of failure cascade ($N_p = 4$) and $3$ disabled assets ($N_c = 3$) and randomly choose two assets ($f=25$ and $f=30$) to plot the differences $P_{s,f,i}^* - P_{s,f,i-1}^*$ for $i \in \{2,3,4\}$ in Figures~\ref{fig:weights_change_1} and~\ref{fig:weights_change_2}. We label the values of $P_{s,f,i}^* - P_{s,f,i-1}^*$ by ``Stage\_i-1\_i'' in the legend of both figures, and each figure reports the differences for one asset. The horizontal axis in both figures shows all the upstream assets on which the chosen asset is dependent. We see that the majority of the differences between the optimal arc weights in two consecutive stages are zero while the maximum absolute difference is only slightly larger than $0.05$. This observation suggests that, although the differences are relatively small,  our setting provides more robustness than does the more classical setting of considering a single worst-case realization of the arc weights.  

\subsubsection{Cascading failures}
\begin{figure*}[t!]
    \centering
    \begin{minipage}[t]{0.5\textwidth}
        \centering
        \includegraphics[width = \linewidth]{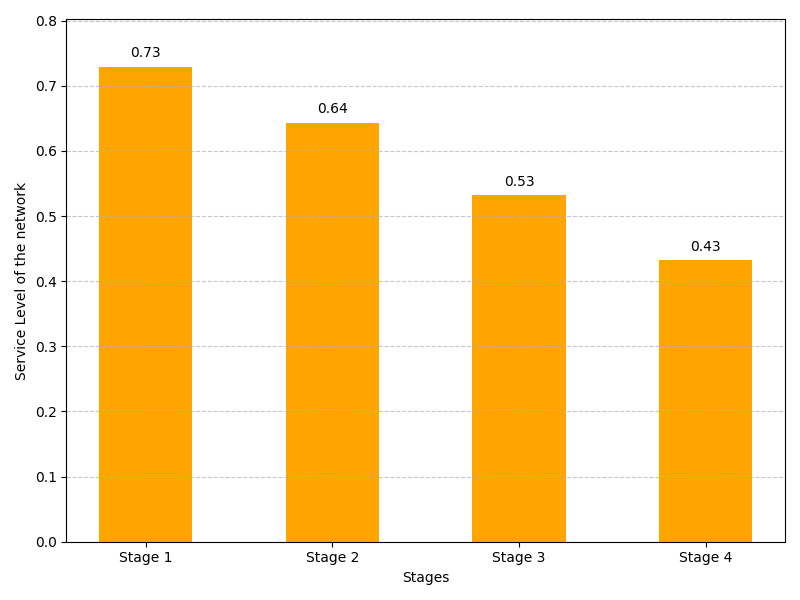}
        \caption{Service level of the infrastructure network in each stage of failure cascade for the network \texttt{a128} with $N_p = 4$ and $N_c = 3$}\label{fig:total_service_level}
    \end{minipage}%
    ~ 
    \begin{minipage}[t]{0.5\textwidth}
        \centering
        \includegraphics[width = \linewidth]{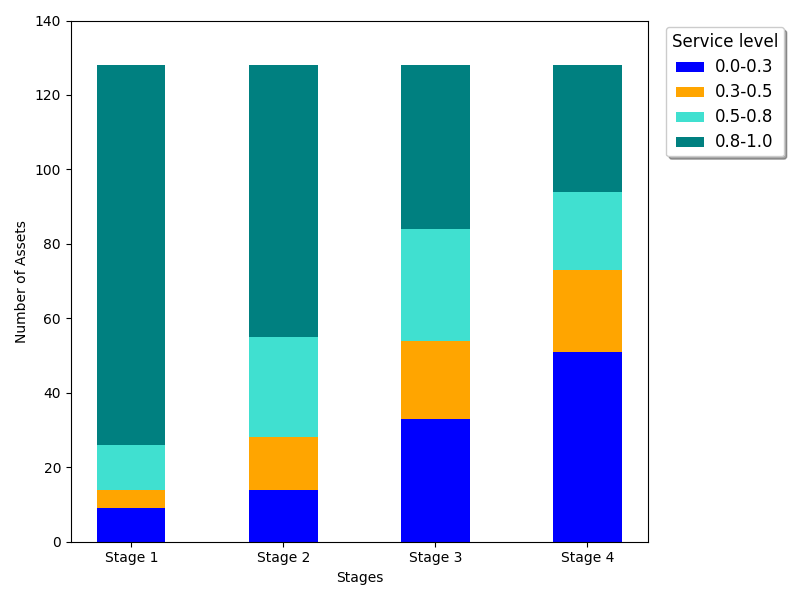}
        \caption{Distribution of service level of individual assets in each stage of failure cascade for the network \texttt{a128} with $N_p = 4$ and $N_c = 3$}\label{fig:service_level_count}
    \end{minipage}
\end{figure*}

We now demonstrate the cascading failure effect as the propagation of failures by comparing the service levels of the infrastructure network as a whole and the distributions of individual service levels in each stage of failure cascade. In particular, we consider the network \texttt{a128} with $4$ stages of failure cascade ($N_p = 4$) and $3$ disabled assets ($N_c = 3$). Figure~\ref{fig:total_service_level} reports the service level of the infrastructure network in each stage of failure cascade. For each stage $i \in \{1,2,3,4\}$, we plot the service level of the network defined by 
\begin{equation}
    \text{service level of the network} = \frac{\sum_{f \in \mathcal{F}} W_f y_{f,i}^*}{\sum_{f \in \mathcal{F}} W_f},
\end{equation}
where $y_{f,i}^*$ represents the optimal value of the variable $y_{f,i}$. We  see a clear trend of decreasing service level of the network as the failures propagate. Additionally, Figure~\ref{fig:service_level_count} presents the distribution of service level of individual assets in each stage of failure cascade, thus also demonstrating the effect of cascading failures. The majority of the assets have a service level larger than $0.8$ in stage 1; but cascading failures lead to degradation of service levels, and a large portion of the assets have a service level less than $0.3$ in stage 4. 

\subsection{Comparison with existing work}
As we discussed in Section \ref{sec:lit_review},  two categories of approaches exist. One category relies on reformulation and relaxation of the original problem to obtain a single-level problem. As we consider the worst-case scenario from the network, the resulting follower problem is a bilinear multistage problem. Even with the equivalent formulation $(P_{mc})$ with the McCormick envelopes, it is difficult to reformulate the bilevel problem into a single-level problem. The size of the final problem could also be too large for good performances. The other approach involves decomposition and iterative procedures. We first discuss the iterative algorithm proposed in \cite{Tang_2016} and \cite{Taninmis_2021} that leverages an approximate leader problem constructed from previously explored leader problem solutions. We observed in our numerical study of this algorithm that it has to explore every single feasible leader problem solution before it terminates. Because of the large number of feasible leader problem solutions in our setting, the algorithm cannot obtain the optimal solution with the target gap of $1\%$ even for the network \texttt{a128} with one stage of failure cascade ($N_p = 1$) in the given time limit. The main procedure in the improved version of the algorithm is to construct a new follower problem solution from a sequence of previously explored solutions. The new follower problem solution has to be feasible to all constraints, and thus a valid lower bound can be obtained. However, the construction of such a follower problem solution assumes certain structures on the follower problem, and our follower problem does not share such structures. Similarly, the Benders algorithm developed in \cite{israeli.wood:02} focuses on the network-interdiction problem of maximizing the shortest path, whereas our problem cannot be readily posed for the algorithm to be applicable. Next, we note that \texttt{Benders} may share some similarity with the algorithm proposed in \cite{Salmeron_2009} and \cite{Salmeron_2015}; nonetheless, the strengthened variant \texttt{Benders-S} has showed further improvements in its performance. 


\section{Conclusion and future work} \label{sec:conclusion}
We study the effect of cascading failures in an interdependent infrastructure a system of interdependent infrastructure systems. In particular, we try to identify the most severe contingencies in the network in the event of a natural disaster that could disable some of the assets and cause failures to propagate further to other dependent assets. We propose two equivalent bilevel interdiction models for this problem, taking into account  unknown or stochastic factors affecting the interconnectivity of the assets with an uncertainty set on the arc weights. 
Because of the complexity of the follower problem, it is difficult to reformulate the bilevel problem into a single-level problem. We instead propose a Benders-type decomposition algorithm and study a strengthened variant of it. We validate the proposed algorithm using two infrastructure networks generated from anonymized real-world data and vary the number of stages of failure cascades and number of initially disabled assets. The proposed algorithm is shown to be able to obtain the optimal solutions in most cases and feasible solutions of practical usefulness in some challenging cases. 

We mention a few future directions for this research that remain to be explored. One direction is exploring other ways to strengthen the optimality cuts and study other classes of valid inequalities to help close the gaps for the challenging cases when the number of stages of failure cascades is large. Furthermore, we can study and take advantage of the structures in the network and dependency graph to further reduce the feasible space in the leader problem. One example is to utilize any symmetries in the network and dependency graph to generate more valid inequalities in each iteration of the algorithm.

\noindent \textbf{Acknowledgement}
\small{
This material is based upon work supported by the 
U.S. Department of Energy,
Office of Science, 
Office of Advanced Scientific Computing Research, 
under Contract DE-AC02-06CH11357. This work was also supported by
the U.S. Department of Energy through grant DE-FG02-05ER25694.}

\bibliographystyle{siam} 
\bibliography{CascadingFailure} 

\noindent \textbf{Government license}

\small{
\noindent The submitted manuscript has been created by UChicago Argonne, LLC, Operator of Argonne National Laboratory ("Argonne”). Argonne, a U.S. Department of Energy Office of Science laboratory, is operated under Contract No. DE-AC02-06CH11357. The U.S. Government retains for itself, and others acting on its behalf, a paid-up nonexclusive, irrevocable worldwide license in said article to reproduce, prepare derivative works, distribute copies to the public, and perform publicly and display publicly, by or on behalf of the Government. The Department of Energy will provide public access to these results of federally sponsored research in accordance with the DOE Public Access Plan (http://energy.gov/downloads/doe-public-access-plan).}


\end{document}